\documentclass{notices}
\usepackage{amsfonts,amssymb,amsmath,amscd,graphicx,amsrefs}
\def\real{{\mathbb R}}
\def\complex{{\mathbb C}}
\def\B{{\cal B}}
\def\sinc{\hbox{\rm sinc}}
\def\qed{\hbox{\vrule width 2.5pt depth 2.5 pt height 3.5 pt}}

\title{Unbounded Growth of Band-Limited Functions}

\author{
Lloyd N. Trefethen
\affil{Professor of Applied Mathematics in Residence,
School of Engineering and Applied Sciences, Harvard University}
}

\begin{document}

\maketitle

Let $f$ be a real or complex function of $\real$.  The idea of $f$
being {\em band-limited} is that it is composed of components
$e^{ikx}$ spanning a finite range of wave numbers $k$.  To be
concrete, let us take the range to be $k\in [-1,1]$.  The basic
example of a function with wave numbers in this range is the
{\em sinc function,}
\begin{equation}
\sinc(x) = {\sin(x)\over x},
\end{equation}
with $\sinc(0) = 1$ and $\sinc(n\pi) = 0$ for the other integers $n$,
as illustrated in Figure~1.

\begin{figure}[h]
\vspace*{15pt}
\noindent\kern 1pt\includegraphics [clip, width=2.9in]{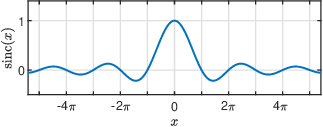}
\vspace*{-5pt}
\caption{The sinc function is the basic example of a band-limited function, and
the starting point of sampling theory.}
\end{figure}

A fundamental result for making such ideas precise is the
{\em Paley-Wiener theorem\/} for band-limited
functions in $L^2$~\cite[Thm.~11.1.2]{simon}, \cite{havin}.
If $F\in L^2(\real)$ has
compact support in $[-1,1]$, then its inverse Fourier transform
\begin{equation}
f(x) = \int_{-\infty}^\infty F(k) \kern 1pt e^{ikx} dk
\label{ft}
\end{equation}
belongs to $L^2(\real)$ and extends
to an entire function of $z = x+iy$ satisfying
\begin{equation}
|f(z)| \le C e^{|y|}, \quad z\in\complex\kern 1pt 
\label{bound}
\end{equation}
for some $C$.
(Throughout this note, $C$ is a generic positive constant, changing
from one appearance to the next.)
Conversely, if $f\in L^2(\real)$
extends to an entire function of $z = x+iy$ that satisfies
(\ref{bound}) for some $C$, then its Fourier transform
\begin{equation}
F(k) = {1\over 2\pi} \int_{-\infty}^\infty f(x) \kern 1pt e^{-ikx} dx
\label{ift}
\end{equation}
belongs to $L^2(\real)$ and has compact support
in $[-1,1]$.  In the case of
$\sinc(x)$, the Fourier transform is $\pi$ times the characteristic
function of $[-1,1]$, and since $\sin(x) = (e^{ix}-e^{-ix})/2\kern
.3pt i$, it is obvious that $\sinc(x)$ satisfies (\ref{bound}).

Band-limited functions became a subject of concerted
attention especially through the work of Henry Landau,
Henry Pollak, and David Slepian at Bell Labs in the early 1960\kern .4pt s.
Like their colleagues Claude Shannon and Richard Hamming, these
men were concerned with
{\em sampling theory,} the study of relationships beteween a continuous
signal and its samples taken at regular intervals.  The {\em Nyquist sampling
theorem} asserts that an $L^2$ function $f$ that is
band-limited to $[-1,1]$ can be recovered
from its samples $\{f(\pi n)\}$, and this observation highlights
the importance of the sinc function: the recovery formula is
\begin{equation}
f(x) = \sum_{n=-\infty}^\infty f(\pi n) \kern 1.5pt \sinc(x-\pi n).
\label{recovery}
\end{equation}

This brings us to a puzzle I encountered
two or three years ago, which led to a conjecture.  The application
area is {\em numerical analytic continuation\/} of an analytic
function $f$ beyond the real or complex domain where it is known.
Analytic continuation is a standard notion that can
be effected in theory by the Weierstrass chain-of-disks method,
but what about practical algorithms applicable to functions just
known numerically?  The best general method seems to be to make
use of rational approximations to $f$, and from this work emerged
a curious empirical observation.

Here is the {\em one-wavelength principle}~\cite{jjiam}.  It seems
that in the usual $16$-digit computer arithmetic, all kinds of
oscillatory functions can be numerically analytically continued about
one wavelength beyond their domain of definition (assuming this domain is
big enough) before accuracy is lost.  This is a rough observation,
not tied to any precise notions of accuracy or wavelength.
The appearance of the number {\em one\/} 
in the principle is a coincidence
related to computations typically employing approximations to 13
digits of accuracy: with 26 or 39 digits one can track about 2
or 3 wavelengths, respectively.  Figure 2 illustrates the effect.
In this example, the function $f(x)$ being extrapolated is a portion
of the trajectory of one of the three components of a solution to
the Lorenz equations.
More examples are given in~\cite{jjiam}.

\begin{figure}[h]
\vspace*{15pt}
\begin{center}
\noindent\kern 1pt\includegraphics [clip, width=2.8in]{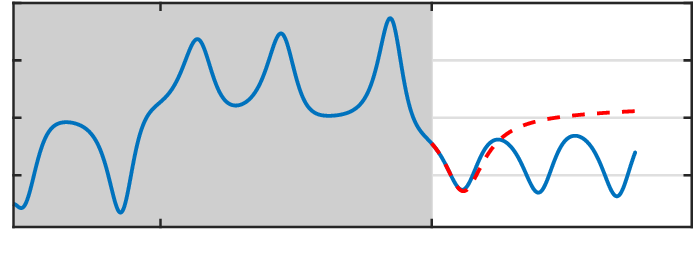}
\end{center}
\vspace*{-15pt}
\caption{The blue curve is an analytic function $f(x)$.
Numerical analytic continuation~\cite{jjiam} to the right of
the gray interval gives the dashed red curve, which
matches $f(x)$ for about one wavelength.}
\end{figure}

The one-wavelength principle is a rough
observation, and I tried to develop a possible
explanation of it by investigating
what seemed a simple model of a problem of this kind.
Let $\B$ be the class of
functions $f\in L^2(\real)$
that are band-limited to $[-1,1]$ and satisfy $|f(x)|\le 1$ for $x\le 0$. 
Define
\begin{equation}
M(x) = \sup_{f\in \B} |f(x)|, \quad x > 0.
\label{Mdef}
\end{equation}
How fast does $M$ grow as a function of $x$?  
The conjecture was that it grows exponentially as $x\to\infty$.
For suppose this were true with an approximate growth rate
\begin{equation}
M(x)\approx e^{Cx}, \quad
C \approx {\log(10^{13})\over 2\pi} \approx 4.8,
\label{rate}
\end{equation}
which corresponds to a factor of
$10^{13}$ over an interval of length $2\pi$.
Then this would give some kind of explanation of the observed behavior.

A 1986 paper by Landau~\cite{landau} (1931--2020) is the sole publication I
have found that considers $M(x)$, and it is encouragingly
consistent with the conjecture.  Landau states the following
theorem in the language of sampling theory.
\begin{quotation}
\noindent {\em Theorem 2:} When sample measurements are accurate
only to within $\varepsilon > 0$ in amplitude or in total energy,
good extrapolation is possible for only a bounded distance (having
an order of magnitude $-\log \varepsilon$) beyond the interval of
observation, regardless of the amount of data used.
\end{quotation}
Landau does not give quantitative bounds, but his proof of the theorem
can be unwound to show that $M(x)$ grows at least exponentially
with a constant $C$ about $1/8$ of that proposed in (\ref{rate}).
He established this exponential growth by considering translates of
the band-limited functions $\sinc(x), (\sinc(x/2))^2, (\sinc(x/3))^3,
\dots.$

Since 1986, Landau's paper has had impact in discussions of
{\em superresolution,} {\em superoscillation,} and {\em prolate
spheroidal wave functions,} which had been introduced by Landau, Pollak, and
Slepian themselves in their earlier
work; a few of the many references in
this area are~\cite{berry,roadmap,lindberg,rokhlin,rogers}.
Landau's theorem is cited for example as establishing that extrapolations of
band-limited functions can be
``exponentially unstable''~\cite{lindberg} or 
``highly unstable''~\cite{rogers}.  All these areas of discourse are
concerned with matters of how rapidly a band-limited function can change
from one region to another, putting them very much in the realm of the
$e^{Cx}$ conjecture.

To try to find a proof and a derivation of $C$, I discussed the
problem with various colleagues.  One of these was
John Urschel at MIT, who shared it with Alex Cohen,
a graduate student at MIT working with Larry Guth.\ \ Cohen,
to my surprise, saw that
the conjecture is false.  In fact, there is no bound on how fast
band-limited functions can grow.  Here is Cohen's result.

\vskip 5pt
{\em {\bf Theorem.}  With the definition $(\ref{Mdef})$,
$M(x) = \infty$ for all $x>0$.}
\vskip 5pt

\noindent In a word, if you can bound a band-limited function
on one side, that gives no constraints on its magnitude.  (It's
different if the function is bounded on both sides.)  Thus the
one-wavelength principle, if it is true in some sense, must find
another mathematical explanation.

\smallskip

{\em Cohen's proof.} 
Consider the function
\begin{equation}
g(x) = \cos(a\sqrt{-x}\kern 1.5pt) = \cosh(a\sqrt x\kern 1.5pt),
\label{ga}
\end{equation}
where $a>0$ is a parameter.  This function is entire (the evenness of
$\cosh$ takes care of the square root) and it satisfies $|g(x)|
\le 1$ for $x \le 0$.  For any $x>0$, it grows without bound
as $a\to\infty$.
What's missing is that $g$ is not in $L^2(\real)$, so it is not
in the class $\B$.  However, this can be fixed by 
multiplying it by a rapidly-decaying band-limited function $\psi$, so that
in fact our counterexample becomes the $a$-dependent family of
functions
\begin{equation}
f(x) = g(x) \psi(x).
\label{prod}
\end{equation}
Specifically, we choose
$\psi$ to be a nonzero entire function satisfying
\begin{equation}
|\psi(z)| \le e^{|y|/2 - |x|^\sigma}, \quad z= x+iy \in\complex\kern 1pt 
\label{bound2}
\end{equation}
for some $\sigma\in(1/2,1)$.
Since $g(x) = O(\exp(a|x|^{1/2}))$ and $\sigma > 1/2$,
such a choice guarantees that $f\in L^2(\real)$ for all
$a$.
The reason for requiring $\sigma< 1$ is that this allows
such a bandlimited function to exist.  There can be no 
band-limited function with decay $\psi(x) = O(\exp(-C|x|))$, because
its Fourier transform would have to
be analytic (by another Paley-Wiener theorem),
which precludes compact support.
But it is known that band-limited functions exist satisfying
(\ref{bound2}) for any $\sigma < 1$~\cite{beurling}, \cite[Thm.~1.4.1]{bjorck}.
To show that $M(x) = \infty$ for a given $x>0$ and thus
prove the theorem, we just
have to pick such a $\psi$ that is nonzero at this value of
$x$ to ensure that $|f(x)|\to\infty$ as $a\to\infty$.

To finish the argument, it remains to confirm that $f$
is band-limited as
required, satisfying the condition (\ref{bound}). 
Following (\ref{ga}), (\ref{prod}), and (\ref{bound2}), we calculate
\begin{displaymath}
|f(x+iy)| \le \exp( a |x+iy|^{1/2} - |x|^\sigma + |y|/2) .
\label{bound3}
\end{displaymath}
Now since
$|x+iy|\le 2\max\{|x|,|y|\}$ and therefore
$|x+iy|^{1/2}\le 2(|x|^{1/2}+|y|^{1/2})$,
we have
\begin{displaymath}
|f(x+iy)| \le \exp( 2a |x|^{1/2} + 2a |y|^{1/2}- |x|^\sigma + |y|/2) .
\end{displaymath}
The proof is completed by noting that for some constant $C$,
\begin{displaymath}
\exp( 2a |x|^{1/2} + 2a |y|^{1/2}- |x|^\sigma ) 
\le C \exp(|y|/2).
\end{displaymath}
The value of $C$ depends on $a$, but uniformity with respect to
$a$ is not needed.  \qed
\medskip

The function $\cosh(a\sqrt x \kern 1pt)$ has remarkable properties.
For $x\le 0$ it is just an oscillatory cosine of an argument
varying with $x$, but the oscillation gets faster as $a$ increases, as
shown in Figure~3.  One would hardly guess from the plots that this whole $a$-dependent
class of functions is uniformly band-limited, but this is the case---and
precisely so in the $L^2$ sense once $g$ is multiplied by $\psi$.  What
makes this possible is the exponentially
great scale of $g(x)$ for $x>0$, growing rapidly
as $a\to \infty$.  Thus $x\le 0$ lies at the edge of the main
signal, and as is well known to experts in superoscillation and
prolate spheroidal wave functions, almost anything is possible in the edges.

\begin{figure}[h]
\vspace*{15pt}
\noindent\kern 6pt\includegraphics [clip, width=2.9in]{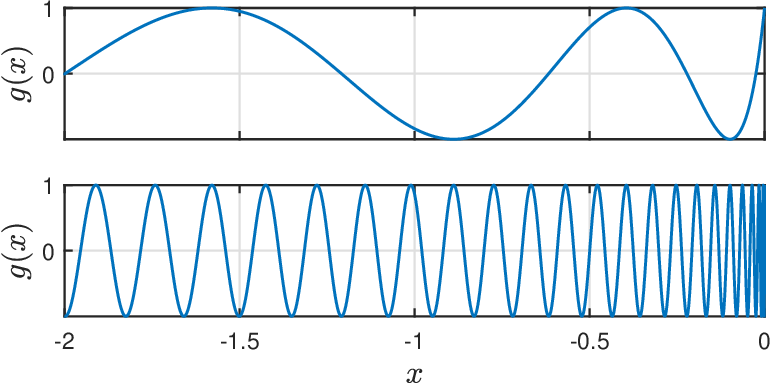}
\vspace*{-5pt}
\caption{The function $g(x)$ of $(\ref{ga})$ on the negative real axis
for $a=10$ (above) and $100$ (below).
Despite the increasingly rapid oscillations as $a\to\infty$, all such functions
are uniformly band-limited after multiplication by the fixed envelope $\psi(x)$,
an example of superoscillation.  For $x>0$, they take huge values.}
\end{figure}

\medskip
{\bf Acknowledgments.}  Of course the main thanks must go to Alex
Cohen for sharing his theorem and the
striking example of $\cosh(a \sqrt x\kern 1pt)$.
In addition I am grateful for advice from
Michael Berry, Greg Beylkin, Robert Calderbank, Karlheinz
Gr\"ochenig, Achim Kempf, Lucas Monz\'on, Vladimir Rokhlin,
Gaurav Thakur, and John Urschel.

\smallskip

\bibliographystyle{plain}
\bibliography{refs}

\end{document}